\documentclass{article}
\usepackage{latexsym}
\usepackage{amssymb}
\usepackage{amsfonts}
\usepackage{amsmath}

\newtheorem{defi}{Definition}[section]
\newtheorem{coro}[defi]{Corollary}

\newtheorem{lema}[defi]{Lemma}
\newtheorem{prop}[defi]{Proposition}
\newtheorem{teor}[defi]{Theorem}
\newtheorem{obs}[defi]{Remark}

\newlength{\soraya}     

\newcommand{\hm}{\hspace*{.5cm}}

\newcommand{\huc}{\hspace*{.1cm}}

\newcommand{\dm}{\noindent {\bf Proof:} }

\newcommand{\su}{\mbox{$\huc\subseteq\huc$}}

\newcommand{\bc}{\begin{center}}
\newcommand{\ec}{\end{center}}

\newcommand{\IF}{{I\kern-0.3emF}}
\newcommand{\IK}{{I\kern-0.3emK}}

\newcommand{\IQ}{{I\kern-0.3emQ}}
\newcommand{\IR}{{I\kern-0.3emR}}

\newcommand{\V}{\mathcal V}

\newcommand{\ia}{\'\i}

\newcommand{\ra}{\rightarrow}

\newcommand{\cpr}{\copyright}

\newcommand{\n}{\neg}

\newcommand{\rar}{\rightarrow}

\newcommand{\srel}{\stackrel}
\newlength{\chico}
\newcommand{\seq}[2]{
\settowidth{\chico}{$\srel{\textstyle \underline{\hspace{.3cm} #1
\hspace{.3cm}}}{#2}$}
\parbox{\chico}{\raisebox{-1cm}{$\srel{\textstyle
\underline{\hspace{.3cm} #1 \hspace{.3cm}}}{#2}$}}\,}

\def\pushright#1{{\parfillskip=0pt\widowpenalty=10000
  \displaywidowpenalty=10000\finalhyphendemerits=0
  \leavevmode\unskip\nobreak\hfil\penalty50\hskip.2em\null
  \hfill{#1}\par}}

\newcommand{\curtains}{\pushright{$\Box$}\par \vspace{1ex}}
\begin{document}
%
\title{ A General Axiomatization for the logics\\
of the Hierarchy ${\mathbb{I}}^n {\mathbb{P}}^k$ \footnote{Researth supported by CICITCA - National University of San Juan.}}
\author{V\ia ctor Fern\'andez\\
\small Basic Sciences Institute (Mathematical Area)\\ 
\small Philosophy College\\
\small National University of San Juan\\
\small Av. Ignacio de la Roza 230 (Oeste) \\
\small San Juan, Argentina \\
\small E-mail:~{\tt vlfernan@ffha.unsj.edu.ar}
}
\date{}
\maketitle

\begin{abstract} In this paper, the logics of the family ${\mathbb{I}}^n {\mathbb{P}}^k$:=$\{{ I^n P^k}\}_{(n,k) \in \omega^2}$ are formally defined by means of finite matrices, as a simultaneous generalization of the weakly-intuitionistic logic $I^1$ and  of the paraconsistent logic $P^1$. It is proved that this family can be naturally ordered, and it is shown an adequate axiomatics for each logic of the form $I^n P^k$. 
\end{abstract}

\



\normalsize

\section {Introduction and preliminaries}

\noindent The propositional logic ${P^1}$ was defined by A. Sette in \cite{set:73}, within the context of a wide research about Paraconsistent Logic developed in the 70's. It possesses special characteristics that distinguish it from the family $\{C^n\}_{\{0 \leq n \leq \omega\}}$, the fundamental paraconsistent hierarchy (see \cite{dac:74}). Among other properties, even when $P^1$ can be defined by means of a Hilbert-Style axiomatics, it can be also obtained by means of a finite matrix (meanwhile no one of the $C^n$-logics can be characterized in this way). The matrix semantics for $P^1$ is built taking as basis a set of three truth-values: $T_0$ and $F_0$ (intended as the ``classical truth-values'')  together with $T_1$ (which can be associated to an ``intermediate truth'').  Besides that, $P^1$ is {\it maximal} w.r.t to the propositional classical logic (${ CL}$), in the sense that, if any axiom-schema (independent of the original ones) is added to the axiomatics of ${P^1}$, then this new axiomatics generates ${CL}$. Finally, $P^1$ is {\it algebraizable}, as it was shown in \cite{lew:mik:sch:90}. 

\

\noindent As a dual counterpart of this logic, A. Sette and W. Carnielli defined in \cite{set:car:95} the logic ${I^1}$, which, in general terms, shares with $P^1$ several properties among the already mentioned (finite axiomatizability, maximality relative to ${CL}$ and algebraizability). Besides that, one of the more remarkable differences between ${I^1}$ and ${P^1}$ is the following: in ${ P^1}$ is not valid the {\it non-contradiction principle} $NCP$: $\n (\n \phi \wedge \phi)$, but it holds the {\it middle excluded principle} $MEP$: $\n \phi \vee \phi$. On the other hand, $I^1$ behaves exactly in the opposite way: it verifies $NCP$ and it does not verify $MEP$.

\

\noindent The logic ${I^1}$ is defined by means of a $3$-valued matrix, too: in this case (and unlike $P^1$), the ``new truth value'' is $F_1$, an ``intermediate truth-value of falsehood''. Considering this fact, it was suggested in \cite{set:car:95} a generalization of these logics {\it by the addition of new intermediate truth-values}, in such a way that the ``new logics'' already obtained constitute a family (which could be ordered in a natural way). Following (and simplifying at some extent) these suggestions, it was defined in \cite{fer:01} the family  ${\mathbb{I}}^n {\mathbb{P}}^k$ = 
$\{ { I^n P^k} \}_{(n,k) \in \omega^2}$. Every member of ${\mathbb{I}}^n {\mathbb{P}}^k$
(usually mentioned here just as an {\it ${ I^n P^k}$-logic}) can be considered as a generalization of  ${I^1}$ and of ${P^1}$ at the same time, by several reasons. First of all, the classical logic $CL$ can be identified simply with ${I^0 P^0}$. Similarly, $P^1$ (resp. $I^1$) is simply
${I^0 P^1}$ (resp. ${I^1 P^0}$). Moreover, every ${ I^n P^k}$-logic has $n+k+2$ truth-values (as it will be seen later). In addition, it can be estabilished an order relation within ${\mathbb{I}}^n {\mathbb{P}}^k$ . The logics of this family fail to verify $NCP$ and/or $MEP$ (with the obvious exception of $I^0 P^0$ that satisfies both properties). It is worth to comment that, since the ${ I^n P^k}$-logics are 
finite-valued, and (mostly) paraconsistent/weakly-intuitionistic ones, they can be used as ``laboratory logics'' in the study of several interesting properties (see  \cite{car:con:mar:07} or \cite{ram:fer:09}, for example). 

\

\noindent However, an open problem referred to this family consists of providing an adequate (i.e. sound and complete) axiomatics for {\it all} the ${ I^n P^k}$-logics. This paper is essentially devoted to offer a suitable axiomatics for them. Moreover, the soundness and completeness theorems shown here can be considered {\it general} in this sense: their proofs are given in such a way that the adequacity of all the logics of ${\mathbb{I}}^n {\mathbb{P}}^k$ (w.r.t. to the axiomatics here presented) is demonstrated in a structured mode, common to any pair $(n,k) \in \omega^2$ previously fixed. The technique to prove this result is adapted to the well-known Kalm\'ar{'}s method to prove completeness for $CL$ (see \cite{men:97}).

\

\noindent To avoid unnecesary information or formalism, the notions to be used to prove adequacity will be reduced as much as possible (this entails that this paper will contain some notational abuses). Besides that, the structure of this article is as follows: in the next section the ${ I^n P^k}$-logics will be defined by means of finite matrices, some simple properties will be shown here, and it will be defined an order relation $\preceq$ in the family 
${\mathbb{I}}^n {\mathbb{P}}^k$ (this justifies the expression ``hierarchy'' used for this family). In addition, it will be demonstrated that ${I^{n_1} P^{k_1}} \preceq 
{I^{n_2} P^{k_2}}$ if and only if $(n_2, k_2) \leq (n_1, k_1)$. In Section \ref{axiomatica-inpk}, it will presented a general axiomatics for all the $I^n P^k$-logics and it will proven some properties, which are essential to the proof of adequacity (result developed in Section \ref{adecuacidad-inpk}). For that, it is assumed that the reader is familiarized with the notions of formal proof, schema axioms, inference rules 
and so on, within the context of Hilbert-Style axiomatics. So, the definitions of these concepts (and other related ones) will be omitted. This paper concludes with some comments about future work. 

\section{Semantic Presentation of the Hierarchy ${\mathbb{I}}^n {\mathbb{P}}^k$} \label{matrices-para-inpk} 

To define a matrix semantics for the logics of the family ${\mathbb{I}}^n {\mathbb{P}}^k$, it is necessary to start with the definition of the language $L(C)$, common to all the $I^n P^k$-logics:

\begin{defi}\label{lenguaje-inpk}\rm{The {\bf set of connectives of all the $I^n P^k$-logics} is $C$:=$\{\n,\rar\}$, with obvious arities. The {\bf language $L(C)$} (or {\it set of formulas}) for the $I^n P^k$-logics is the algebra of words generated by $C$ over a countable set $\V$, in the usual way.}
\end{defi}

\noindent Along this paper, the uppercase greek letters $\Gamma, \Delta, \Sigma \dots$ denote {\it sets of formulas} of $L(C)$. In addition, the lowercase greek letters $\phi$, $\psi$, $\theta$ are metavariables ranging over the {\it individual formulas of $L(C)$}. Finally,  the letters $\alpha$, $\alpha_1$, $\alpha_2, \dots$ will be used as metavariables referred only to the {\it atomic formulas} (that is, the elements of $\V$). All these notations can be used with subscripts, when neccesary. 
On the other hand, the expression $\phi[\alpha_1,\dots,\alpha_m]$ indicates that the atomic formulas ocurring on $\phi$ {\it are precisely} $\alpha_1,\dots,\alpha_m$ (this expression will be applied in the developement  of the completeness proof). 

\

\noindent Despite their common language, the difference between each one of the $I^n P^k$-logics is given by their respective matrix semantics, defined as follows:

\begin{defi}\label{logica-inpk} \rm{Let $(n,k) \in \omega^2$, with $\omega$ = $\{0,1,2,\dots\}$. The 
{\bf matrix $M(n,k)$} is defined as a pair 
$M_{(n,k)}$ = $(({A_{(n,k)}},C_{(n,k)}),
D_{(n,k)})$, where:

\noindent {\bf a)} $(A_{(n,k)},C_{(n,k)})$ is an algebra, {\it similar to $L(C)$}, with support 

\begin{center}$A_{(n,k)}$: = $\{ F_0, F_1, \cdots, F_n, T_0, T_1, \cdots
T_k \}$ \footnote{Every algebra $(A_{(n,k)},C_{(n,k)})$ will be identified with its support, if there is no risk of confussion. }\end{center}

\noindent {\bf b)} $D_{(n,k)}$ = $\{ T_0, T_1, \cdots T_k \}$.

\noindent In addition, The operations $\n$ and $\rar$ of $C_{(n,k)}$ (also called {\bf truth-functions}) \footnote{Strictly speaking, the operations of $C_{(n,k)}$ {\it are not} the connectives of $C$, of course. However, they will be denoted in the same way for the sake of simplicity.} are defined by the truth tables indicated below.}
\end{defi}

\bc $
\begin{array}{|c|c c c c|} \hline
 & F_0 & F_r & T_i & T_0 \\ \hline \neg & T_0 & F_{r-1} & T_{i-1}
 & F_0 \\ \hline
\end{array}
$

\

\

$
\begin{array}{|c|c c c c|} \hline
\rar & F_0 & F_s & T_j & T_0 \\ \hline F_0 & T_0 & T_0 & T_0 & T_0\\
F_r & T_0 & T_0 & T_0 & T_0\\
T_i & F_0 & F_0 & T_0 & T_0\\
T_0 & F_0 & F_0 & T_0 & T_0 \\ \hline
\end{array}
$

\

\

\hspace{0.5 cm} \mbox{\it With $1
\le r,s \le n;\hspace{0.2 cm}1 \le i,j \le k$ .} \ec

\begin{obs}\label{explicacion-tabla-inkp}\rm{
Realize that
the truth-values $F_1,\dots,F_n$ can be considered informally as {\it intermediate values of falsehood}, meanwhile $T_1,\dots,T_k$ are {\it intermediate values of truth}. 
In addition, every application of $\n$ to a ``non classical value'', approximates more and more the value to the ``classical ones'', 
$F_0$ and $T_0$. Note that there are needed $n$ negations at most to pass from $F_r$ to $F_0$. Similarly, the values of the form $T_i$ ``become'' $T_0$ after $k$ negations at most. On the other hand, the implication $\rar$ cannot distinguish between classical or intermediate truth-values: it just considers every value of the form $F_i$ as being $F_0$, and every value of the form $T_j$ as being $T_0$.} 
\end{obs}

\noindent Taking into account the previous truth-tables, some secondary (and useful) truth-functions can be defined. As a motivation, it would be desirable that disjunction ($\vee$) and conjunction ($\wedge$) behave as $\rar$ in this aspect: they cannot distinguish classical from intermediate truth-values. For that, it is taken as starting point the unary function of ``classicalization'' $\copyright$ (the meaning of this neologism is obvious), defined by $\copyright (A)$ := $(A \rar A) \rar A$, for every $A \in A_{(n,k)}$. So, the truth-table associated to it is

\bc $
\begin{array}{|c|c c c c|} \hline
 & T_0 & T_i & F_r & F_0 \\ \hline \copyright & T_0 & T_0 & F_0
 & F_0 \\ \hline
\end{array}
$ \ec

\noindent  From $\copyright$ it is defined the truth-function $\sim$, of {\it strong} (also called {\it classical}) negation, as $\sim A$ : = $\n ( \copyright A)$. So, its associated truth-table is

\bc   $\begin{array}{|c|c c c c|} \hline
 & F_0 & F_r & T_i & T_0 \\ \hline \sim & T_0 & T_0 & F_0
 & F_0 \\ \hline
\end{array}
$ \ec 

\noindent It is possible to define $\vee$ and $\wedge$ now, adapting the usual definition for $CL$:
$A \vee B$:= $\sim A \rar B$, meanwhile
$A \wedge B$:= $\sim (A \rar \, \sim B)$.
For these connectives, their associated truth-functions are:

\bc $
\begin{array}{|c|c c c c|} \hline
\vee & F_0 & F_s & T_j & T_0 \\ \hline F_0 & F_0 & T_0 & T_0 & T_0\\
F_r & F_0 & F_0 & T_0 & T_0\\
T_i & T_0 & T_0 & {T_0} & {T_0}\\
T_0 & T_0 & T_0 & {T_0} & {T_0} \\ \hline
\end{array}
$ \hspace{2cm} 
$
\begin{array}{|c|c c c c|} \hline
\wedge & F_0 & F_s & T_j & T_0 \\ \hline F_0 & F_0 & F_0 & {F_0} & {F_0}\\
F_r & F_0 & F_0 & {F_0} & {F_0}\\
T_i & {F_0} & {F_0} & {T_0} & {T_0}\\
T_0 & {F_0} & {F_0} & {T_0} & {T_0} \\ \hline
\end{array}
$
\ec

\hspace{0.5 cm} \bc \mbox{With $1 \le i,j \le k;\hspace{0.2 cm}
1 \le r,s \le n$.} \ec

\noindent From the previous definitions, it is clear that all the binary truth-functions consider all the non-designated values $F_j$ as behaving  as $F_0$, and similarly for all the values $T_i$. The same fact holds for $\sim$. 
In the case of $\n$, however, all the truth-values can be differentiated. This is the main difference of $\n$ and $\sim$, and justifies the definition and the study of the $I^n P^k$-logics. For example, when $n \geq 1$, $MEP$: $\n \phi \vee \phi$ it is not an $I^n P^k$-tautology (it is enough to consider 
a valuation $v$ such that $v(A)$ = $T_i$ with $i \geq 1$), meanwhile this principle is valid if $\n$ is replaced by $\sim$. 
That is, $\models_{(n,k)} \sim \phi \vee \phi$ for any $I^n P^k$-logic. In a dual way, when $k \geq 1$, $\sim (\sim \phi \, \wedge \phi)$ is a $I^n P^k$-tautology. (for every $(n,k) \in \omega^2$), but $\n (\n \phi \wedge \phi)$ is  not valid in all the $I^n P^k$-logics. Indeed, $\n (\n \phi \wedge \phi)$ is only valid in the  $I^n P^0$-logics.

\

\noindent After a deeper analysis it is possible to see that, given a fixed logic $I^n P^k$, $\models_{(n,k)}\n \phi \vee \phi$ iff  $\phi$ = $\n^t (\phi \rar \theta)$ (with $t \geq 0$), or $\phi$ = $\n^t \alpha$, with $\alpha \in \V$, with $t \geq n$: otherwise (when $\phi$ = $\n^t \alpha$ with $t < n$) $\not\models_{(n,k)}\n \phi \vee \phi$.  In a similar way, $\models_{(n,k)} \n (\n \phi \wedge \phi)$ iff $\phi$ = $\n^t (\psi \rar \theta)$ with $t \geq 0$, or $\phi$ = $\n^t \alpha$, with $t \geq k$, $\alpha \in \V$.  From these comments can see that $NCP$ and $MEP$ are not valid in general terms. So, it is natural to distinguish between ``well-behaved'' formulas and ``not well-behaved'' ones (with respect to each of the mentioned principles). This distinction is formalized with the unary ``well-behavior'' truth-functions, defined in the obvious way: $A^{\ast}$:= $\n A \vee A$; $A^{\circ}$:= $\n (\n A \wedge A)$, for every $A \in A_{(n,k)}$. Its respective truth-tables are

\bc $
\begin{array}{|c|c c c c|} \hline
 & F_0 & F_r & T_i & T_0 \\ \hline ^{\ast} & T_0 & F_0 & T_0
 & T_0 \\ \hline
\end{array}
$  \hspace{2cm}$
\begin{array}{|c|c c c c|} \hline
 & F_0 & F_r & T_i & T_0 \\ \hline
 ^{\circ} & T_0 & T_0 & F_0
 & T_0 \\ \hline
\end{array}
$ \ec

\

\noindent Besides the behavior of the truth-function in each matrix $M_{(n,k)}$, recall that its definition is motivated by the definition of {\it a consequence relation on $L(C)$} (and therefore of a {\it logic}), in the usual way:

\begin{defi} \label{cons-matricial-inpk}\rm{
An {\bf $M_{(n,k)}$-valuation} is any homomorphism $v:L(C) \to A_{(n,k)}$ (this notion makes sense because $L(C)$ and ${A_{(n,k)}}$ are similar algebras). Recall here that every $M_{(n,k)}$-valuation can be defined just considering functions $v:\V \to A_{(n,k)}$ and extending it to all $L(C)$. The {\bf logic $I^n P^k$} is the pair $I^n P^k$:=$(C,\models_{(n,k)})$, being $\models_{(n,k)} \subseteq \wp(L(C)) \times L(C)$ defined as usual: $\Gamma \models_{(n,k)}\phi$ iff, for very $M_{(n,k)}$-valuation $v$, $v(\Gamma) \subseteq D_{(n,k)}$ implies $v(\phi) \in D_{(n,k)}$. In this context, $\phi$ is an {\bf $I^n P^k$-tautology} iff $\emptyset \models_{(n,k)} \phi$ (this fact will be denoted by $\models_{(n,k)}\phi$, as usual). 
The family $\{I^n P^k\}_{(n,k) \in \omega^2}$ 
will be denoted by ${\mathbb{I}}^n {\mathbb{P}}^k$.
}
\end{defi}

\begin{obs}\rm{The family ${\mathbb{I}}^n {\mathbb{P}}^k$. includes some well-known logics. Indeed, $I^0 P^0$ is just the classical logic  $CL$. On the other hand, the logic $I^1 P^0$ is $I^1$ indeed meanwhile $I^0 P^1$ is just $P^1$. In addition, all the $I^n P^k$-logics can be ``naturally ordered'', taking into account the following definition:
}\end{obs}

\begin{defi}\label{jerarquia-inpk}{\rm The {\bf order relation} $\preceq \su ({\mathbb{I}}^n {\mathbb{P}}^k)^2$ is defined 
in the following natural way: $I^{n_1} P^{k_1}\preceq I^{n_2} P^{k_2}$ iff, for every $\Gamma \cup \{\phi\} \su L(C)$, $\Gamma \models_{(n_1,k_1)}\phi$ implies $\Gamma \models_{(n_2,k_2)}\phi$.}
\end{defi}

 \noindent Taking into account the previous definition, it is natural to visualize $({\mathbb{I}}^n {\mathbb{P}}^k, \preceq)$ as a lattice:
 
\begin{prop}\label{inpk-clave}\rm{In the logic $I^n P^k$ ($n$, $k$ fixed), the following formulas are tautologies:

\noindent a) $\n^{n+1}\phi \vee \n^n \phi$ \hspace{2cm}($(n+1)$-generalization of $MEP$.)

\noindent b) $\n (\n^{k+1}\phi \wedge \n^k{\phi})$ \hspace{1.3cm} ($(k+1)$-generalization of $NCP$).}
\end{prop}
%
%
%
%
%
%
%
\begin{prop}\rm{$I^{n_1}P^{k_1}\preceq I^{n_2}P^{k_2}$ iff $(n_2, k_2)  \leq_{\Pi} (n_1,k_1)$ (being $\leq_{\Pi}$ the order of the product on $\omega^2$). Therefore, the {\it Hierarchy $({\mathbb{I}}^n {\mathbb{P}}^k, \preceq)$} is a lattice.
}
\end{prop}

\noindent \dm  If $(n_2,k_2) \leq_{\Pi} (n_1,k_1)$, then $A_{(n_2,k_2)} 
\su
A_{(n_1,k_1)}$, and $D_{(n_2,k_2)} \su D_{(n_1,k_1)}$. Now suppose that $\Gamma_0 \not\models_{(n_2,k_2)}\phi_0$ for some $\Gamma_0 \cup \{\phi_0\} \su L(C)$. So, there exists a valuation $v:\V \to A_{(n_2,k_2)}$ such that $v(\Gamma_0) \su D_{(n_2,k_2)}$, $v(\phi_0) \notin D_{(n_2,k_2)}$. Define the valuation $w:\V \to A_{(n_1,k_1)}$ as $w(\alpha)$ = $v(\alpha)$, for every $\alpha \in \V$. It can be proved that, for every $\psi \in L(C)$, $w(\psi)$ = $v(\psi)$. Thus, $w(\Gamma_0) \su D_{(n_2,k_2)} \su D_{(n_1,k_1)}$ and $w(\phi_0) \in \{F_0,\dots,F_{n_2}\} \su \{F_0,\dots,F_{n_1}\}$. That is, 
$\Gamma_0 \not\models_{(n_1,k_1)}\phi_0$. The previous argument shows that $I^{n_1} P^{k_1}\preceq I^{n_2} P^{k_2}$.

For the converse, suppose $(n_2,k_2) \not\leq_{\Pi} (n_1,k_1)$. There are two cases that must be analyzed in different ways. 
First, if $n_2 > n_1$ consider any formula $\phi_1$:= $\n^{n_1 + 1} \alpha \vee \n^{n_1}\alpha$, with $\alpha \in \V$. So, $\models_{(n_1,k_1)}\phi_1$, by Prop. \ref{inpk-clave}.a). Now, defining the valuation $v_1:\V \to A_{(n_2,k_2)}$ by $v_1(\alpha)$:=$F_{n_2}$, it holds $v_1(\phi_1)$ = 
$\n^{n_1 +1} F_{n_2} \vee \n^{n_1} F_{n_2}$ = $F_{n_2 - (n_1 + 1)} \vee F_{n_2 - n_1}$ = $F_0$ (since $n_1 + 1 \leq n_2$). Thus, $\not\models_{(n_2,k_2)}\phi_1$. On the other hand, if $k_2 > k_1$, let $\phi_2$:=$\n (\n^{k_1 + 1} \alpha \wedge \n^{k_1}\alpha)$. As in the first case,  
$\models_{(n_1,k_1)}\phi_2$, by Prop. \ref{inpk-clave}.b). Now, if it is defined the valuation $v_2:\V \to A_{(n_2,k_2)}$ such that 
$v_2(\alpha)$ = $T_{k_2}$, then $\not\models_{(n_1,k_1)}\phi_2$ (note here that $k_1 + 1 \leq k_2$). So, for both possibilities
it holds $I^{n_1} P^{k_1} \not\preceq I^{n_2} P^{k_2}$. This concludes the proof.\curtains

\noindent Some consequences of the previous result, useful to visualize $\preceq$ (actually, its underlying strict order $\prec$) are the following:

\begin{coro}\rm{In ${ I^n P^k}$ it holds that:

\noindent a)  $I^{n+1} P^{k} \prec I^n P^k$.

\noindent b)  $I^n P^{k+1} \prec I^n P^k$.

\noindent c) $I^n P^{k+1}$ and $I^{n+1} P^{k}$ are not comparables.
}
\end{coro}

\noindent  This section concludes with the mention of the following result that will be applied at the end of this paper:

\begin{prop}\label{models-finitaria-dt}\rm{The consequence relation $\models_{(n,k)}$ verifies:}

\noindent a) $\Gamma \models_{(n,k)} \phi$ implies $\Gamma \cup \{\psi\} \models_{(n,k)} \phi$ \hfill [Monotonicity]

\noindent b) $\Gamma, \phi \models_{(n,k)}\psi$ iff $\Gamma \models_{(n,k)}\phi \rar \psi$ \hfill[Semantic Deduction Theorem]

\noindent c) If $\Gamma \models_{(n,k)}\phi$, then $\Gamma{'}\models_{(n,k)}\phi$ for some finite set $\Gamma{'}\su \Gamma$ \hfill[Finitariness]
 
\end{prop}

\dm Obviously, it is holds a). The claim b) arises from the truth-table of $\rar$. With respect to c), $\models_{(n,k)}$ is finitary because is naturally defined by means of a single finite matrix (result indebted to R. W\'ojcicki: see \cite{woj:88}). \hfill $\Box$

\section{A Hilbert-Style Axiomatics for the $I^n P^k$-logics}\label{axiomatica-inpk}

From now on, consider an $I^n P^k$-logic fixed, with $(n,k) \in \omega^2$. To obtain the desired axiomatics, the secondary truth-functions $\sim$, $\copyright$, $\vee$ and $\wedge$ from the previous section will be reflected by means of the definition of secondary connectives in $L(C)$. 
Formally:

\begin{defi} \label{conectivos-secundarios}\rm{The secondary connectives $\sim$, $\copyright$, $\vee$, $\wedge$ are defined in $L(C)$ in the following way:

\noindent $\copyright \phi$:= $(\phi \rar \phi) \rar \phi$

\noindent $\sim \phi$:= $\n (\copyright \phi)$

\noindent $\phi \vee \psi$:= $\sim \phi \rar \psi$

\noindent $\phi \wedge \psi$:= $\sim(\phi \rar \sim \psi)$. 

\noindent $\phi^{\ast}$:= $\n \phi \vee \phi$ 

\noindent $\phi^{\circ}$:= $\n( \n \phi \wedge \phi)$ 

\noindent In addition, the conncectives $\vee_{CL}$ and $\wedge_{CL}$ are defined by:

\noindent $\phi \vee_{CL} \psi$:= $\n \phi \rar \psi$

\noindent $\phi \wedge_{CL} \psi$:= $\n(\phi \rar \n \psi)$ \footnote{The ``classical'' connectives $\wedge_{CL}$ and $\vee_{CL}$ are not essential in the proof of Completeness. However, they are indicated here for a better explanation of the comparison between these connectives w.r.t $\wedge$ and $\vee$, as it will be remarked later.}. 

\noindent Finally, the expression $\n^q \phi$ indicates $\n( \dots (\n \phi))\dots)$, $q$ times. $\n^0 \phi$ is merely $\phi$. 
}
\end{defi}

\noindent Taking into account the previous conventions, the axiomatics for the $I^n P^k$-logics will be presented in the sequel. For that consider, from now on, an arbitrary 
(fixed) pair $(n,k)\in \omega^2$.

\begin{defi}\label{axiomatica-inpk-definitiva}\rm{The {\bf consequence relation} $\vdash_{(n,k)} \su \wp(L(C))\times L(C)$ is defined by means of the following Hilbert-Style axiomatics, considering these schema axioms:

\

\noindent $Ax_1$ $\phi \rar (\psi \rar \phi)$

\noindent $Ax_2$ $(\phi \rar (\psi \rar \theta)) \rar ((\phi \rar \psi) \rar (\phi \rar \theta))$

\noindent $Ax_3$ $(\phi \rar \psi)^{\ast}$

\noindent $Ax_4$ $(\phi \rar \psi)^{\circ}$

\noindent $Ax_5$ $(\n^n \phi)^{\ast}$

\noindent $Ax_6$ $(\n^k \phi)^{\circ}$

\noindent $Ax_7$ $\phi^{\ast} \rar [\psi^{\circ} \rar ((\n \phi \rar \n \psi) \rar ((\n \phi \rar \psi) \rar \phi))]$

\noindent $Ax_8$ $\phi^{\ast} \rar [\psi^{\circ} \rar ((\phi \rar \n \psi) \rar ((\phi \rar \psi) \rar \n \phi))]$



\noindent $Ax_9$ $\phi^{\ast} \rar (\n \n \phi \rar \phi)$

\noindent $Ax_{10}$ $\phi^{\circ} \rar (\phi \rar \n \n \phi)$ 

\noindent $Ax_{11}$ $\phi^{\ast}\rar (\n \phi)^{\ast}$ 

\noindent $Ax_{12}$ $\phi^{\circ}\rar (\n \phi)^{\circ}$ 

\

\noindent The only inference rule for this axiomatics is 

\noindent $\mbox{Modus Ponens (MP):}\hm\seq{\phi, \hspace{0.5cm}\phi \ra
\psi}{\psi}.\hspace{0.5cm}$

}
\end{defi}

\

%
%
%

\noindent From this definition, the well-known notions of formal proof (with or without hypotheses), formal theorem, etc. are the usual. 
Because of this, $\vdash_{(n,k)}$ is {\it monotonic}: $\Gamma \vdash_{(n,k)}\phi$ implies $\Gamma \cup \{\psi\} \vdash_{(n,k)}\phi$. This fact will be widely used.

\begin{obs}\label{reflexividad}\rm{ 
It is well known that the inclusion of $Ax_1$, $Ax_2$ and $MP$ entail that it is valid $\vdash_{(n,k)}\phi \rar \phi$. Moreover:}
\end{obs}

\begin{teor} \label{teor-deduccion-sintactico}\rm{$\vdash_{(n,k)}$ satisfies the (syntactic) Deduction Theorem (DT). That is, $\Gamma, \phi \vdash_{(n,k)} \psi$ iff $\Gamma \vdash_{(n,k)}\phi \rar \psi$.}
\end{teor}

\dm This result holds because the inclusion of axioms $Ax_1$ and $Ax_2$ too, and considering that the only (primitive) inference rule is Modus Ponens. See \cite{men:97} for a detailed proof. \hfill $\Box$

\

\noindent $Ax_1$ and $Ax_2$ allow to obtain some useful rules in relation to $\vdash_{(n,k)}$, too:

\begin{prop}\label{reglas-secundarias}\rm{Given the logic $I^n P^k$, the following secondary rules are valid:}
\end{prop}

\noindent $\mbox{Permutation (Perm):}\hm\seq{\phi \rar (\psi \rar \theta)}{\psi \rar (\phi \rar \theta)}$.

\

\noindent $\mbox{Transitivity (Trans):}\hm\seq{\phi \rar \psi, \hspace{0.5 cm} \psi \rar \theta}{\phi \rar \theta}$.

\

\noindent $\mbox {Reduction (Red):}\hm\seq{(\phi \rar \psi) \rar \theta}{\psi \rar \theta}$.

%
%
%

\

\

\noindent 
The following two results involve formulas of the form $\phi^{\ast}$ or $\phi^{\circ}$:

\begin{prop}\label{simil-idempotencia}\rm{For every $\phi \in L(C)$, for every $(n,k) \in \omega^2$, it holds:

\noindent $\vdash_{(n,k)} (\phi^{\ast})^{\ast}$; $\vdash_{(n,k)} (\phi^{\ast})^{\circ}$; $\vdash_{(n,k)}(\cpr \, \phi)^{\ast}$; $\vdash_{(n,k)} (\cpr \, \phi)^{\circ}$. }
\end{prop}

\noindent This result is valid since $\phi^{\ast}$:= $\sim (\n \phi) \rar \phi$ and $\cpr \, \phi$ = $(\phi \rar \phi) \rar \phi$, and considering axioms $Ax_3$ and $Ax_4$ from Definition \ref{axiomatica-inpk-definitiva}.

\begin{prop} \label{dem-bien-comportadas} \rm{If $\models_{(n,k)}\phi$, then $\vdash_{(n,k)}\phi^{\ast}$ and $\vdash_{(n,k)}\phi^{\circ}$.}
\end{prop}

\noindent \dm If $\models_{(n,k)}\phi$ then (checking the truth-tables of $I^n P^k$) 
$\phi$ is necessarily of the form $\n^q (\psi \rar \theta)$,  with $q \geq 0$.
From this, apply $Ax_3$, $Ax_4$ (and, eventually, $Ax_{11}$ and $Ax_{12}$). \hfill $\Box$

\

\noindent The next result shows some basic {\it $I^n P^k$-theorems}:

\begin{prop} \label{formulas-derivadas-inpk}
\rm{The following formulas of $L(C)$ are 
theorems w.r.t. $\vdash_{(n,k)}$
:}
\end{prop}

\noindent {a)} $\phi \rar \copyright \phi$; \hspace{2cm} {a{'})} $\copyright \phi \rar \phi$ 

\noindent {b)} $\phi^{\ast} \rar (\sim \phi \rar \n \phi)$ 

%

\noindent {c)} $\phi^{\ast}\rar[\psi^{\circ}\rar((\n \phi \rar \n \psi)\rar(\psi \rar \phi))]$

\noindent {d)} $\phi^{\ast}\rar[\psi^{\circ}\rar((\phi \rar \psi)\rar(\n \psi \rar \n \phi))]$

\noindent {e)} $(\sim \phi \rar \sim \psi) \rar ((\sim \phi \rar \copyright \psi)\rar \copyright \phi)$
%
%

\noindent {f)} $(\sim \phi \rar \sim \psi) \rar ((\sim \phi \rar \psi) \rar \phi))$
%
%

\

\dm The following are schematic formal proofs (in the context of $\vdash_{(n,k)}$) for every formula above indicated. Sometimes it will be applied Theorem \ref{teor-deduccion-sintactico} or Proposition \ref{reglas-secundarias} without explicit mention.

\

\noindent For {a)}: $\phi \rar \copyright \phi$ = $\phi \rar ((\phi \rar \phi) \rar \phi)$ is a particular case of $Ax_1$. For the  case of {a{'})}:

\noindent $1)$ $(\phi \rar \phi) \rar  \phi$ \hfill\mbox{[Hyp.; Def. $\copyright \phi$]}\\
\noindent $2)$  $\phi \rar \phi$ \hfill\mbox{[Remark \ref{reflexividad}]}\\
\noindent $3)$ $\phi$  \hfill\mbox{[1), 2), MP]}

\noindent So, it is valid $\copyright \phi \vdash_{(n,k)} \phi$.

\

\noindent For {b)}: 

\noindent $1)$ $\phi^{\ast}$ \hfill\mbox{[Hyp.]}\\
\noindent  $2)$  $(\copyright \phi)^{\circ}$ \hfill\mbox{[Prop. \ref{simil-idempotencia}]}\\
\noindent  $3)$  $\phi^{\ast} \rar [\, (\copyright \phi)^{\circ} \rar ((\phi \rar \sim \phi) \rar ((\phi \rar 
\copyright \phi) \rar \n \phi)]$ \hfill\mbox{[$Ax_8$, 
Def. \ref{conectivos-secundarios} (of $\sim$)]}\\
\noindent  $4)$  $(\phi \rar \sim \phi) \rar ((\phi \rar \copyright \phi)\rar \n \phi)$ \hfill\mbox{[$1)$, $2)$, $3)$, MP 
]}\\
\noindent  $5)$   $(\phi \rar \copyright \phi) \rar ((\phi \rar \sim \phi)\rar \n \phi)$ \hfill\mbox{[$4)$, Perm.]}\\
\noindent  $6)$   $\phi \rar \copyright \phi$ \hfill\mbox{[{a)}]}\\
\noindent  $7)$   $(\phi \rar \sim \phi) \rar \n \phi$ \hfill\mbox{[$5)$, $6)$, MP]}\\
\noindent  $8)$   $\sim \phi \rar (\phi \rar \sim \phi)$ \hfill\mbox{[$Ax_1$
]}\\
\noindent  $9)$  $\sim \phi \rar \n \phi$ \hfill\mbox{[$7)$, $8)$, Trans.]}

\noindent That is, $\phi^{\ast} \vdash_{(n,k)} \sim \phi \rar \n \phi$. 

\

%
%
%
%
%
%

\noindent For {c)}:

\noindent $1)$  $\phi^{\ast}$ \hfill\mbox{[Hyp.]}\\
\noindent $2)$  $\psi^{\circ}$ \hfill\mbox{[Hyp.]}\\
\noindent $3)$  $\n \phi \rar \n \psi$ \hfill\mbox{[Hyp.]}\\
\noindent $4)$ $\psi$ \hfill\mbox{[Hyp.]}\\
\noindent $5)$ $\psi \rar (\n \phi \rar \psi)$ \hfill\mbox{[$Ax_1$]}\\
\noindent $6)$ $\n \phi \rar \psi$ \hfill\mbox{[4), 5), MP]}\\
\noindent $7)$ $\phi^{\ast}\rar[\psi^{\circ}\rar((\n \phi \rar \n \psi)\rar ((\n \phi \rar \psi) \rar \phi))]$ \hfill\mbox{[$Ax_7$]}\\
\noindent $8)$ $(\n \phi \rar \n \psi)\rar ((\n \phi \rar \psi) \rar \phi)$ \hfill\mbox{[7), 1), 2), MP]}\\
\noindent $9)$ $(\n \phi \rar \psi) \rar \phi$ \hfill\mbox{[8), 3), MP]}\\
\noindent $10)$ $\phi$  \hfill\mbox{[9), 6), MP]}

\noindent Thus, $\phi^{\ast},\psi^{\circ},\n \phi \rar \n \psi, \psi \vdash_{(n,k)}\phi$. 

\

\noindent For {d)}: adapt the proof of c), using $Ax_8$ instead of $Ax_7$. Then, it will be valid
%
%
$\phi^{\ast},\psi^{\circ},\phi \rar \psi, \n \psi \vdash_{(n,k)}\n \phi$. 

\

\noindent For {e)}:

\noindent $1)$  $(\copyright \phi)^{\ast}$ \hfill{[Prop. \ref{simil-idempotencia} ]}\\
\noindent $2)$  $(\copyright \psi)^{\circ}$ \hfill{[Prop. \ref{simil-idempotencia}]}\\
\noindent $3)$   $(\copyright \phi)^{\ast} \rar [(\copyright \psi)^{\circ} \rar ((\sim \phi \rar \sim \psi) \rar ((\sim \phi \rar \copyright \psi) \rar \copyright \phi))]$ 

\hfill{[Def. \ref{conectivos-secundarios} (of $\sim$), $Ax _7$.]}\\
\noindent $4)$  $(\sim \phi \rar \sim \psi) \rar ((\sim \phi \rar \copyright \psi) \rar \copyright \phi)$ \hfill{[1), 2), 3), MP]}

\noindent So, $\vdash_{(n,k)} (\sim \phi \rar \sim \psi) \rar ((\sim \phi \rar \copyright \psi) \rar \copyright \phi)$. 
\
%
%
%
%
%

\

\noindent For {f)}:

\noindent $1)$  $\sim \phi \rar \sim \psi$ \hfill\mbox{[Hyp.]}\\
\noindent $2)$  $\sim \phi \rar \psi$ \hfill\mbox{[Hyp.]}\\
\noindent $3)$  $\psi \rar \copyright \psi$ \hfill\mbox{[{a)}]}\\
\noindent $4)$  $\sim \phi \rar \copyright \psi$ \hfill\mbox{[2), 3), Trans.]}\\
\noindent $5)$  $(\sim \phi \rar \sim \psi) \rar ((\sim \phi \rar \copyright \psi) \rar \copyright \phi))$ \hfill\mbox{[{e)}]}\\
\noindent $6)$  $\copyright \phi$ \hfill\mbox{[1), 4), 5), MP]}\\
\noindent $7)$  $\copyright \phi \rar \phi$ \hfill\mbox{[{a{'})}]}\\
\noindent $8)$  $\phi$ \hfill\mbox{[6), 7), MP]}

\noindent From all this, $\sim \phi \rar \sim \psi, \sim \phi \rar \psi \vdash_{(n,k)}\psi$. Now, apply Theorem \ref{teor-deduccion-sintactico}, as in the previous results. This concludes the proof. \hfill $\Box$

\begin{obs} \label{vee-wedge-inpk} \rm{Now is convenient to relate the axiomatics given in Definition \ref{axiomatica-inpk-definitiva} with a well-known axiomatics for $CL$ = $I^0 P^0$. According \cite{men:97}, $CL$ can be axiomatized by MP joined with the following three schema axioms:

\noindent $Bx_1$ = $Ax_1$

\noindent $Bx_2$ = $Ax_2$

\noindent $Bx_3$ = $(\n \phi \rar \n \psi) \rar ((\n \phi \rar \psi) \rar \phi)$. 

\noindent Note that, cf. Definition \ref{axiomatica-inpk-definitiva}, fixed an arbitrary consequence relation $\vdash_{(n,k)}$, the axiom $Bx_3$ of the previous axiomatics is replaced by  a weaker version ($Ax_7$). Anyway, since in the particular case of $\vdash_{(0,0)}$, axioms $Ax_5$ and $Ax_6$ establish that, for every formula $\phi \in L(C)$, $\vdash_{(0,0)}\phi^{\ast}$ and 
$\vdash_{(0,0)}\phi^{\circ}$, it is possible to recover the axiomatics determined by $Bx_1$, $Bx_2$ and $Bx_3$, actually. Moreover:
}
\end{obs}

\begin{prop}\label{sim-comporta-como-n}\rm{Let $\phi \in L(C)$, in such a way that $\phi$ is a formal theorem of $CL$ (that is, $\vdash_{(0,0)} \phi$), and let $\phi{'} \in L(C)$, obtained by $\phi$ replacing all the ocurrences of the symbol $\n$ in $\phi$ by $\sim$. Then $\vdash_{(n,k)}\phi{'}$.}
\end{prop}

\dm Consider the axiomatics for $I^0 P^0$ indicated in Remark \ref{vee-wedge-inpk}, and compare it with the general axiomatics given in Definition \ref{axiomatica-inpk-definitiva}. First of all note that neither $Bx_1$ (= $Ax_1$) nor $Bx_2$ (= $Ax_2$) have ocurrences of $\n$. Besides, 
since $Bx_3$=$(\n \phi \rar \n \psi) \rar ((\n \phi \rar \psi) \rar \phi))$, and considering Prop. \ref{formulas-derivadas-inpk}.f), it holds
$\vdash_{(n,k)} (\sim \phi \rar \sim \psi) \rar ((\sim \phi \rar \psi) \rar \phi)$ ( = ${Bx_3}{'}$). From these facts, it can be easily proved by induction of the lenght of the formal proof of $\phi$ (w.r.t $\vdash_{(0,0)}$)  that $\vdash_{(0,0)}\phi$ implies $\vdash_{(n,k)}\phi{'}$. \hfill $\Box$

%
\begin{coro}\label{sim-comporta-como-n-y-0}\rm{Suppose $\phi \in L(C)$, and
let the formula $\phi{''} \in L(C)$ built by $\phi$ replacing the eventual ocurrences of $\n$ in $\phi$ by $\sim$, and replacing  every ocurrence of $\vee_{CPL}$ (resp. $\wedge_{CPL}$), understood as an abbreviation (cf. Definition \ref{conectivos-secundarios}), by $\vee$ (resp. 
$\wedge$). Then, $\vdash_{(0,0)}\phi$ implies $\vdash_{(n,k)}\phi{''}$. 
}
\end{coro}

\noindent For instance, since $\vdash_{(0,0)} \n \phi \vee_{CL} \phi$, then $\vdash_{(n,k)} \sim \phi \vee \phi$. However, it is not generally valid that 
$\vdash_{(n,k)} \n \phi \vee_{CL} \phi$, obviously. The next result collects some particular cases of the previous corollary:

\begin{coro}\label{ejemplos-sim-n-y-o}\rm{The relation $\vdash_{(n,k)}$ verifies, given $(n,k) \in \omega^2$:

%
\noindent {a)} 
$ \vdash_{(n,k)} \phi  \rar \phi \vee \psi$; \hspace{2.1cm} {a{'})} $\vdash_{(n,k)} \psi \rar \phi \vee \psi$

\noindent {b)} 
$ \vdash_{(n,k)} \phi \wedge \psi \rar \phi$; \hspace{2.1cm} {b{'})} $\vdash_{(n,k)} \phi \wedge \psi \rar \psi$ 

\noindent {c)} 
$ \vdash_{(n,k)} (\phi \rar \theta) \rar ((\psi \rar \theta) \rar (\phi \vee \psi \rar \theta))$

\noindent {d)}  
$\vdash_{(n,k)}\phi \rar (\psi \rar (\phi \wedge \psi))$.

\noindent {e)}  
$\vdash_{(n,k)}\phi \wedge \psi \rar \phi \vee \psi$.
}
\end{coro}

\noindent Finally, to prove Completeness, it will be necessary:

\begin{prop}\label{lista-formulas-utiles-2}\rm{The following are $I^n P^k$-theorems:}
\end{prop}

\noindent a) $\phi \rar \phi^{\ast}$ 

\noindent b) $\phi^{\circ} \rar ( \n \phi \rar (\phi \rar \psi))$

\noindent c) $ (\phi^{\circ})^{\circ}$

\noindent d) $\n (\phi^{\ast})\rar \phi^{\circ}$ 
%

\noindent e) $\, \sim \phi \rar \phi^{\circ}$

\noindent f) $\phi^{\ast} \rar (\n (\phi \vee \psi) \rar \n \phi)$
%

\noindent g) $\psi^{\circ} \rar ( \phi \rar (\n \psi \rar \n (\phi \rar \psi)))$

\noindent h) $(\n(\phi \rar \psi))^{\ast}$; \hspace{0.5cm} $ (\n(\phi \rar \psi))^{\circ}$

\noindent i) $(\n (\phi^{\ast}))^{\circ}$

\noindent j) $\n (\phi^{\ast}) \rar (\phi \rar \psi)$
%
%

\

\dm these formal theorems are formally demonstrated as in Proposition \ref{formulas-derivadas-inpk}, applying DT and Proposition \ref{reglas-secundarias} if were necessary: 

\

\noindent For a): taking into account Corollary \ref{ejemplos-sim-n-y-o}.a{'}), it is valid $\vdash_{(n,k)} \phi\rar \n \phi \vee \phi$. That  is, $\vdash_{(n,k)} \phi \rar \phi^{\ast}$.

\

\noindent For b):

\noindent $1)$  $\phi^{\circ}$ \hfill\mbox{[Hyp.]}\\
\noindent $2)$  $\n \phi$ \hfill\mbox{[Hyp.]}\\
\noindent $3)$  $(\n \phi \rar \psi)^{\ast}$ \hfill\mbox{[$Ax_3$]}\\
\noindent $4)$  $(\n \phi \rar \psi)^{\ast}\rar[\phi^{\circ} \rar((\n(\n \phi \rar \psi) \rar \n \phi)\rar(\phi \rar (\n \phi \rar \psi)))]$ 

\hfill\mbox{[Prop. \ref{formulas-derivadas-inpk}.c)]}\\
\noindent $5)$  $(\n(\n \phi \rar \psi) \rar \n \phi)\rar(\phi \rar (\n \phi \rar \psi))$ \hfill\mbox{[1), 3), 4), MP]}\\
\noindent $6)$  $\n \phi \rar (\n(\n \phi \rar \psi) \rar \n \phi)$ \hfill\mbox{[$Ax_1$]}\\
\noindent $7)$  $\n(\n \phi \rar \psi) \rar \n \phi$   \hfill\mbox{[2), 6), MP]}\\
\noindent $8)$  $\phi \rar (\n \phi \rar \psi)$ \hfill\mbox{[5), 7),MP]}\\
\noindent $9)$  $\n \phi \rar (\phi \rar \psi)$ \hfill\mbox{[8), Perm.]}\\
\noindent $10)$ $\phi \rar \psi$ \hfill\mbox{[2), 9), MP]}

\noindent Thus, it holds $\phi^{\circ} \vdash_{(n,k)} \n \phi \rar (\phi \rar \psi)$, as was desired.

\

\noindent For c): 

\noindent $1)$  $(\copyright (\n \phi \rar \sim \phi))^{\circ}$ \hfill\mbox{[Prop. \ref{simil-idempotencia}]}\\
\noindent $2)$  $(\sim (\n \phi \rar \sim \phi))^{\circ}$ \hfill\mbox{[$1)$, $Ax_{12}$]}\\
\noindent $3)$  $(\n \phi \wedge \phi)^{\circ}$ \hfill\mbox{[$2)$, Def. \ref{conectivos-secundarios} (of $\wedge$)]}\\
\noindent $4)$  $(\n (\n \phi \wedge \phi))^{\circ}$ \hfill\mbox{[$3)$, $Ax_{12}$]}\\
\noindent $5)$  $(\phi^{\circ})^{\circ}$ \hfill\mbox{[$4)$, Def. \ref{conectivos-secundarios} (of $\circ$)]}

\

\noindent For d): 

\noindent $1)$ $(\phi^{\ast})^{\circ}$ \hfill\mbox{[Prop. \ref{simil-idempotencia}]}\\
\noindent $2)$ $(\copyright(\n \phi \rar \sim \phi))^{\ast}$ \hfill\mbox{[Prop. \ref{simil-idempotencia}]}\\
\noindent $3)$ $(\n \phi \wedge \phi)^{\ast}$ \hfill\mbox{[2), Def. \ref{conectivos-secundarios} (of $\wedge$), $Ax_{11}$]}\\
\noindent $4)$ $(\n \phi \wedge \phi)^{\ast} \rar [(\phi^{\ast})^{\circ} \rar((\n \phi \wedge \phi \rar \phi^{\ast}) \rar ((\n (\phi^{\ast}) \rar \phi^{\circ}))]$
\hfill\mbox{[Prop. \ref{formulas-derivadas-inpk}.{d)}]}\\
\noindent $5)$ $(\n \phi \wedge \phi \rar \n \phi \vee \phi) \rar ((\n (\phi^{\ast}) \rar \phi^{\circ})$ \hfill\mbox{[1), 3), 4), MP]}\\
\noindent $6)$ $\n \phi \wedge \phi \rar \n \phi \vee \phi$ \hfill\mbox{[Corollary \ref{ejemplos-sim-n-y-o}.e)]}\\
\noindent $7)$ $\n (\phi^{\ast}) \rar \phi^{\circ}$ \hfill\mbox{[5), 6), MP]}

\noindent So, $\vdash_{(n,k)} \n (\phi^{\ast}) \rar \phi^{\circ}$

\

%
%
%
\noindent For e):

\noindent $1)$ $(\n \phi \wedge \phi)^{\ast} \rar [(\copyright \phi)^{\circ} \rar ((\n \phi \wedge \phi \rar \copyright \phi) \rar (\sim \phi \rar \phi^{\circ}))]$

\hfill\mbox{[Prop. \ref{formulas-derivadas-inpk}.d), Def. \ref{conectivos-secundarios} (of $\sim$ and $\circ$)]}\\
\noindent $2)$ $\n \phi \wedge \phi \rar \phi$ \hfill\mbox{[Corollary \ref{ejemplos-sim-n-y-o}.b{'})]}\\
\noindent $3)$  $\phi \rar \copyright \phi$ \hfill\mbox{[Prop. \ref{formulas-derivadas-inpk}.a)]}\\
\noindent $4)$  $\n \phi \wedge \phi \rar \copyright \phi$ \hfill\mbox{[2), 3), Trans.]}\\
\noindent $5)$  $(\copyright \phi)^{\circ}$ \hfill\mbox{[Prop. \ref{simil-idempotencia}]}\\
\noindent $6)$  $(\copyright(\n \phi \rar \sim \phi))^{\ast}$ \hfill\mbox{[Prop. \ref{simil-idempotencia}]}\\
\noindent $7)$  $(\sim (\n \phi \rar \sim \phi))^{\ast}$  \hfill\mbox{[6), Def. \ref{conectivos-secundarios} (of $\sim$), $Ax_{11}$]}\\
\noindent $8)$ $(\n \phi \wedge \phi)^{\ast}$ \hfill\mbox{[7), Def. \ref{conectivos-secundarios} (of $\wedge$)]}\\
\noindent $9)$  $\sim \phi \rar \phi^{\circ}$ \hfill\mbox{[1), 4), 5), 8), MP]}

\noindent Therefore, $\vdash_{(n,k)}\sim \phi \rar \phi^{\circ}$

\

\noindent For f):

\noindent $1)$  $\phi^{\ast}$ \hfill\mbox{[Hyp.]}\\
\noindent $2)$  $(\phi \vee \psi)^{\circ}$ \hfill\mbox{[$Ax_4$, Def. \ref{conectivos-secundarios} (of $\vee$)]}\\
\noindent $3)$  $\phi^{\ast}\rar[(\phi \vee \psi)^{\circ} \rar((\phi \rar \phi \vee \psi)\rar(\n (\phi \vee \psi)\rar \n \phi)]$ 
\hfill\mbox{[Prop. \ref{formulas-derivadas-inpk}.{d)}]}\\
\noindent $4)$  $(\phi \rar \phi \vee \psi) \rar (\n(\phi \vee \psi) \rar \n \phi)$ \hfill\mbox{[1), 2), 3), MP]}\\
\noindent $5)$  $\phi \rar \phi \vee \psi$ \hfill\mbox{[Corollary \ref{ejemplos-sim-n-y-o}.a)]}\\
\noindent $6)$   $\n(\phi \vee \psi) \rar \n \phi$   \hfill\mbox{[4), 5), MP]}

\noindent Thus, $\phi^{\ast}\vdash_{(n,k)} \n (\phi \vee \psi)  \rar \n \phi$.

%
%
%

\

\noindent For g):

\noindent $1)$  $\psi^{\circ}$ \hfill\mbox{[Hyp.]}\\
\noindent $2)$  $(\phi \rar \psi)^{\ast}$ \hfill\mbox{[$Ax_3$]}\\
\noindent $3)$  $(\phi \rar \psi)^{\ast} \rar [\psi^{\circ} \rar(((\phi \rar \psi) \rar \psi) \rar (\n \psi \rar \n(\phi \rar \psi)))]$ \hfill\mbox{[Prop. \ref{formulas-derivadas-inpk}.d)]}\\
\noindent $4)$  $((\phi \rar \psi) \rar \psi) \rar (\n \psi \rar \n (\phi \rar \psi))$ \hfill\mbox{[1), 2), 3), MP]}\\
\noindent $5)$  $(\phi \rar \psi) \rar (\phi \rar \psi)$ \hfill\mbox{[Remark \ref{reflexividad}]}\\
\noindent $6)$  $\phi \rar ((\phi \rar \psi) \rar \psi)$ \hfill\mbox{[5), Perm.]}\\
\noindent $7)$   $\phi \rar (\n \psi \rar \ (\phi \rar \psi))$   \hfill\mbox{[4), 6), Trans.]}

\noindent So, $\psi^{\circ}\vdash_{(n,k)}\phi \rar (\n \psi \rar \n (\phi \rar \psi))$ is obtained.

\

\noindent For h): By $Ax_3$ and $Ax_{11}$ it holds $\vdash_{(n,k)} (\n(\phi \rar \psi))^{\ast}$; By $Ax_4$ and $Ax_{12}$ it holds $\vdash_{(n,k)} (\n(\phi \rar \psi))^{\circ}$.

\

\noindent For i): It is a particular case of h). 
%
%
%
%

\

\noindent For j):

\noindent $1)$  $\n (\phi^{\ast})$ \hfill\mbox{[Hyp.]}\\
\noindent $2)$  $\phi$ \hfill\mbox{[Hyp.]}\\
\noindent $3)$  $\phi \rar \phi^{\ast}$ \hfill\mbox{[a)]}\\
\noindent $4)$  $\phi^{\ast}$ \hfill\mbox{[2), 3), MP]}\\
\noindent $5)$  $(\phi^{\ast})^{\circ}$ \hfill\mbox{[Prop. \ref{simil-idempotencia} ]}\\
\noindent $6)$  $(\phi^{\ast})^{\circ} \rar (\n(\phi^{\ast}) \rar (\phi^{\ast} \rar \psi))$ \hfill\mbox{[b)]}\\
\noindent $7)$   $\psi$   \hfill\mbox{[1), 4), 5). 6), MP]}

\noindent That is, it holds $\n (\phi^{\ast}), \phi \vdash_{(n,k)}\psi$. Then, apply Theorem \ref{teor-deduccion-sintactico}. This last result completes the proof. \hfill $\Box$

\section {General Soundness and Completeness}\label{adecuacidad-inpk}

It is easy to check that the axioms given in Definition \ref{axiomatica-inpk-definitiva} are $I^n P^k$-tautologies. 
So, taking into account that MP preserves $I^n P^k$-tautologies, it holds:

\begin{teor}\label{inpk-correccion-debil}\rm{[Weak Soundness] If $\vdash_{(n,k)}\phi$, then $\models_{(n,k)}\phi$.}
\end{teor}

\noindent A theorem of (weak) Completeness arises as an adaptation of the well-known Kalm\'ar's proof for Classical Logic $CL$, cf. \cite{men:97}: 

\begin{defi} \label{conjuntos-asociados-inpk}
\rm{For every formula $\phi[{\alpha}_1,{\alpha}_2,\dots,{\alpha}_m] \in L(C)$, for every 
$I^n P^k$-valuation $v$, 
for every atomic formula ${\alpha}_p$  ($1 \leq p \leq m$) let $Q_p^v$ be the set associated to $\alpha_p$ and to $v$, defined  by:
\begin{itemize}
\item If $v({\alpha}_p)$ = $F_r$ (with $1 \leq r \leq n$), then:

$Q_p^v$ =
$\{ \n ({{\alpha}_p}^{\ast}), \n((\n {\alpha}_p)^{\ast}), \dots, \n((\n^{r-1} {\alpha}_p)^{\ast}), (\n^r
{\alpha}_p)^{\ast}\}$
\item If $v({\alpha}_p)$ = $T_i$ (with $1 \leq i \leq k$), then:

$Q_p^v$ =
$\{ \n {\alpha}_p \wedge {\alpha}_p, \n^2 {\alpha}_p \wedge \n {\alpha}_p,\dots,\n^{i}{\alpha}_p \wedge \n^{i-1} {\alpha}_p, (\n^i
{\alpha}_p)^{\circ}\}$
\item If $v({\alpha}_p)$ = $F_0$, then $Q_p^v$= $\{\sim {\alpha}_p, ({\alpha}_p)^{\ast}\}$ 
\item If $v({\alpha}_p)$ = $T_0$, then $Q_p^v$ = $\{\copyright {\alpha}_p, ({\alpha}_p)^{\circ} \}$
\end{itemize}

\noindent In addition, let the set $\Delta_{\phi}^{v}$:= $Q_1^v \cup
Q_2^v,\dots \cup Q_m^v$.

\

\noindent On the other hand, for every $I^n P^k$-valuation $v$ indicated above, the {\bf formula $\phi^v$} (determined by $\phi$ and $v$) is defined as follows:

\begin{itemize}
\item If $v(\phi)$ = $F_r$ (with $0 \leq r \leq n$), then $\phi^v$ = $\n^{r+1} \phi$.
\item If $v(\phi)$ = $T_i$ (with $0 \le i \le k$), then $\phi^v$ = $\phi$.
\end{itemize}
}
\end{defi}

\noindent For the next technical (and essential) result, the following obvious fact will be applied without explicit mention: according to the previous definition, if $\phi \in L(C)$ and 
$\psi$ is a subformula of $\phi$ then, for every valuation $v$, $\Delta_{\psi}^v \su \Delta_{\phi}^v$. Bearing this in mind it is possible to demonstrate:

\begin{lema} \label{tecnico-1} \rm{For every formula $\phi$=$\phi[\alpha_1,\dots,\alpha_m] \in L(C)$, for every $I^n P^k$-va\-lua\-tion $v$, it holds that 
$\Delta_{\phi}^{v} \vdash_{(n,k)} {\phi}^v$.}
\end{lema}

\dm By induction on the complexity of $\phi$. The analysis is divided in the following cases:

\

\noindent \underline{Case 1):} $\phi \in \V$ (without losing generality, $\phi$ = $\alpha_1$, which implies $\Delta_{\phi}^v$ = $Q_1^v$). Then: 

\noindent If $v(\phi)$ = $F_0$, then $\phi^v$ = $\n \phi$ and $\Delta_{\phi}^v$ = $\{\sim \phi, {\phi}^{\ast}\}$. So, $\Delta_{\phi}^v \vdash_{(n,k)}\phi^v$ 
by Prop. \ref{formulas-derivadas-inpk}.b). 

\noindent If $v(\phi)$ = $F_r$ ($1 \leq r \leq n$), then $\phi^v$ = $\n^{r+1} \phi$ and 
$\{(\n^r \phi)^{\ast},
\n (\n^r \phi \vee \n^{r-1} \phi)\} \su \Delta_{\phi}^v$. Now, by Prop. \ref{lista-formulas-utiles-2}.f),  $\vdash_{(n,k)}{(\n^r \phi)}^{\ast}  \rar (\n (\n^r \phi \vee \n^{r-1} \phi) \rar 
\n^{r+1} \phi)$. From all this, $\Delta_{\phi}^v \vdash_{(n,k)} \n^{r+1} \phi$ (= $\phi^v$).

\noindent If $v(\phi)$ = $T_i$ ($1 \leq i \leq k$), then $\n \phi \wedge \phi \in \Delta_{\phi}^v$. By Corollary \ref{ejemplos-sim-n-y-o}.{b{'})}, 
$\Delta_{\phi}^v \vdash_{(n,k)} \phi$ (=$\phi^t$).

\noindent If $v(\phi)$ = $T_0$, then 
$\Delta_{\phi}^t$ = $\{\copyright{\phi},(\phi)^{\circ}\}$. Thus, $\Delta_{\phi}^t \vdash_{(n,k)} \copyright \phi$
( = $(\phi \rar \phi) \rar \phi$)). So, since $\vdash_{(n,k)} \phi \rar \phi$, it holds $\Delta_{\phi}^t \vdash_{(n,k)} \phi$ (= $\phi^t$).
%
%
 The proof of Case 1) is completed.

\

\noindent \underline{Case 2):} when $\phi$ is of the form $\n \psi$. Consider
the following subcases:

\noindent{2.1)}: $v(\psi)$ = $F_0$. By (I.H), $\Delta_{\psi}^v \vdash_{(n,k)} \psi^{v}$ (= $\n \psi$ = $\phi^{v}$). Hence, $\Delta_{\phi}^v \vdash_{(n,k)}\phi^v$.

\noindent{2.2)}: $v(\psi)$ = $F_r$, with $1
\leq r \leq n$. Note that $v(\phi)$ = $F_{r-1}$, which implies $\phi^v$
= $\n^r \phi$ 
= $\n^{r+1} \psi$ = $\psi^v$. So, 
$\Delta_{\phi} \vdash_{(n,k)} \phi^v$, by (I.H).

\noindent{2.3)}:  $v(\psi)$ = $T_1$. 
From Definition \ref{logica-inpk}, $\psi$ = $\n^{q-1}\alpha$, with $1 \leq q \leq k$, $\alpha \in \V$, and $v(\alpha)$ = $T_q$. Thus,
$\Delta_{\psi}^v$ = $Q_{\alpha}^v$ = $\{(\n \alpha \wedge \alpha),\dots,(\n^{q}\alpha \wedge \n^{q-1} \alpha), (\n^q \alpha)^{\circ}\}$. 
From this, 
using Corollary \ref{ejemplos-sim-n-y-o}.b), it holds $\Delta_{\phi}^v \vdash_{(n,k)}\n^q \alpha$ (=$\phi^v$).

\noindent{2.4)}: $v(\psi)$ = $T_i$, with $2
\leq i \leq k$. So, $v(\phi)$ = $T_{i-1}$, $1 \leq i-1 \leq k-1$. 
In addition, $\psi$ = $\n^q \alpha$, with $0 \leq q \leq k-i$, $v(\alpha)$ = $T_{i+q}$ and $\alpha \in \V$
From this, $\phi \wedge \psi$ = $\n^{q+1} \alpha \wedge \n^{q} \alpha 
\in Q_{\alpha}^v$ = $\Delta_{\psi}^v$ = $\Delta_{\phi}^v$ (since $q+1 \leq k$). So, applying Corollary \ref{ejemplos-sim-n-y-o}.b), $\Delta_{\phi}^v\vdash_{(n,k)}\phi$ = $\phi^v$.
%

\noindent 2.5): $v(\psi)$ = $T_0$. Then,
$v(\phi)$ = $F_0$. To prove that $\Delta _{\phi}^v \vdash_{(n,k)} \n
\phi$ (= $\n \n \psi$) it suffices to demonstrate
$\Delta_{\psi}^v \vdash_{(n,k)} \psi^{\circ}$ ${(\star)}$.
Indeed, if this fact holds, from $Ax_{10}$, then it would be verified $\Delta_{\psi}^v \vdash_{(n,k)} \psi \ra \n \n \psi$. And,
since it holds $\Delta_{\psi}^v \vdash_{(n,k)} \psi$ (by (I.H)), it will be obtained
$\Delta_{\phi}^v \vdash_{(n,k)} \phi^v$. Now, to prove ${(\star)}$, consider the following possibilities:


\noindent 2.5.1): $\psi$ is of the form ${\n}^q (\theta_1 \rar \theta_2)$ (with $0 \leq q$). Applying $Ax_4$ and (eventually) $Ax_{12}$, it holds $\vdash_{(n,k)} \psi ^{\circ}$, and therefore $\Delta_{\phi}^v \vdash_{(n,k)} \psi^{\circ}$.

\noindent 2.5.2): $\psi$ is of the form ${\n} ^q \alpha$, $\alpha \in \V$, $0 \leq q$. In this case, $Q_{\alpha}^v$ = $\Delta_{\psi}^v$. Consider here the different
possibilities for $v(\alpha)$:

\noindent 2.5.2.1): $v(\alpha)$ = $F_0$.
Then, $\Delta_{\alpha}^v$ = $\{\sim{\alpha},\alpha^{\ast}\}$. By Prop. \ref{lista-formulas-utiles-2}.e), $\Delta_{\psi}^v \vdash_{(n,k)} \alpha^{\circ}$. Then, apply $Ax_{12}$ ($q$ times).

\noindent 2.5.2.2): $v(\alpha)$ = $F_r$ (with $1 \leq r \le n$). Since $\n(\alpha^{\ast}) \in \Delta_{\alpha}^v$, it holds
$\Delta_{\alpha}^v \vdash_{(n,k)} \alpha^{\circ}$, because Prop. \ref{lista-formulas-utiles-2}.d). From this, $\Delta_{\psi}^v \vdash_{(n,k)} \psi^{\circ}$, by $Ax_{12}$ again. 

\noindent 2.5.2.3): $v(\alpha)$ = $T_i$ (with $1 \leq i \le k $). So, $i \le q$ (in fact: if $i > q$, then 
$v(\psi)$ 
= $T_{i-q}$, $i-q \geq 1$, contradicting  $v(\psi)$ = $T_0$). Besides, note that $(\n^i \alpha)^{\circ} \in
\Delta_{\alpha}^v $, which implies $\Delta_{\psi}^v \vdash_{(n,k)}(\n^i \alpha)^{\circ}$. So,  
since $i \le q$, $\Delta_{\psi} \vdash_{(n,k))}
(\psi)^{\circ}$, again by $Ax_{12}$.

\noindent 2.5.2.4): $v(\alpha)$ = $T_0$. Then, $\Delta_{\psi}^{v} \vdash_{(n,k)}
\alpha^{\circ}$, since $\alpha^{\circ} \in Q_{\alpha}^v \su \Delta_{\psi}^v$. Then, applying $Ax_{12}$  one more time, ${\bf (\star)}$ is valid. The proof of Case 2) is concluded.

\

\noindent \underline{Case 3):} $\phi$ is of the form $\psi
\ra \theta$. There exist the following possibilities \footnote {The
first three subcases indicate the possibilities for $v(\phi)$ =
$T_0$. The last two cases correspond to  $v(\phi)$ = $F_0$.}
:\\
\noindent 3.1): $v(\psi)$ = $F_0$ (and so, $\phi^v$ = $\psi \rar \theta$). By (I.H), $\Delta_{\psi}^v \vdash_{(n,k)} \n \psi$ ($\star$). In addition, it can be proved that
$\Delta_{\psi}^v \vdash_{(n,k)} \psi^{\circ}$ ($\star \star$) (such a  proof runs as follows, according the internal structure of $\psi$):

\noindent 3.1.1): $\psi$ = $\alpha \in \V$. So, 
$\Delta_{\psi}^v \vdash_{(n,k)} \sim \psi$. Applying Prop. \ref{lista-formulas-utiles-2}.e), it holds ($\star \star$). 

\noindent 3.1.2): $\psi$ = $\n^q \alpha$, $1 \leq q$, $\alpha \in \V$. Consider the following possibilities for $v(\alpha)$:

\noindent 3.1.2.1): $v(\alpha)$ = $F_0$. Then, by Prop. \ref{lista-formulas-utiles-2}.e), $Q_{\alpha}^v \vdash_{(n,k)} \alpha^{\circ}$ .

\noindent 3.1.2.2): $v(\alpha)$ = $F_r$ ($1 \leq r \leq n$). Then, $q \geq r$ (because  $q < r$ implies $v(\psi)$ = $v(\n^q \alpha)$ = 
$\n^q (v(\alpha))$ = $F_{r-q}
\neq F_0$, which is absurd). Besides that, $\n ((\n^{r-1} \alpha)^{\ast}) \in Q_{\alpha}^v$. Therefore
$Q_{\alpha}^v \vdash_{(n,k)} (\n^{r-1} \alpha)^{\circ}$, because Prop. \ref{lista-formulas-utiles-2}.d). 

\noindent 3.1.2.3): $v(\alpha)$ = $T_i$ ($1 \leq i \leq k$). So, $q \geq i$ (by similar reasons to 3.1.3.2)).
In addition, $(\n^i \alpha)^{\circ} \in 
Q_{\alpha}^v$, and so $Q_{\alpha}^v \vdash_{(n,k)} (\n^i \alpha)^{\circ}$. 

\noindent 3.1.2.4): $v(\alpha)$ = $T_0$. Obviously, $Q_{\alpha}^v \vdash_{(n,k)} (\alpha)^{\circ}$, from Def. \ref{conjuntos-asociados-inpk}.

\noindent Now note that $Ax_{12}$ can be applied in all the subcases 3.1.2.1)-3.1.2.4), in such a way to obtain $\Delta_{\psi}^{v}\vdash_{(n,k)}\psi^{\circ}$, completing the proof $(\star \star)$ for Subcase 3.1.2).

\noindent 3.1.3): $\psi$ = $\n^q (\theta_1 \rar \theta_2)$, with $0 \leq q$. By $Ax_3$, $\vdash_{(n,k)} (\theta_1 \rar \theta_2)^{\circ}$. Now, apply $Ax_{12}$, 
$q$ times.

\noindent So, it was proven $(\star \star)$ for all the possibilities of Subcase 3.1. From this, $(\star)$ and Prop. \ref{lista-formulas-utiles-2}.b), it holds 
$\Delta_{\phi}^v \vdash_{(n,k)}\psi \rar \theta$ (=$\phi^v$).

\noindent 3.2): $v(\psi)$ = $F_r$, $1 \leq r \leq n$. Again, $v(\phi)$ = $T_0$ and so $\phi^{v}$ =
$\psi \ra \theta$. 
Note that, since $v(\psi)$ = $F_r$, $\psi$ = $\n^q \alpha$, with $q \ge 0$, $\alpha \in \V$. Thus,
$v(\alpha)$ = $F_{r+q}$, with $r+q \le n$, which implies
$Q_{\alpha}^v$ = $\{ \n (\alpha^{\ast}), \n ((\n \alpha)^{\ast}),\dots,
\n ((\n^{r+q-1} \alpha)^{\ast}),(\n^{r+q} \alpha)^{\ast} \}$. So, 
$\n ((\n^q \alpha)^{\ast}) \in
Q_{\alpha}^v$, because $r \geq 1$. Thus, $\Delta_{\phi}^v \vdash_{(n,k)} \n (\psi^{\ast})$. From this and Prop. \ref{lista-formulas-utiles-2}.j), $\Delta_{\phi}^v \vdash_{(n,k)} \phi^v$.

\noindent  3.3): $v(\theta)$ = $T_i$,  $0
\leq i \leq k$. So, 
$\phi^v$ = $\psi \ra \theta$ one more time. By (I.H), $\Delta_{\theta}^v
\vdash_{(n,k)} \theta$. Now apply $Ax_{1}$.

\noindent 3.4): $v(\psi)$ = $T_i$, $v(\theta)$ = $F_r$ ($0 \leq i \leq k$, $1 \leq r \leq n$). 
Using (I.H), $\Delta_{\phi} \vdash_{(n,k)} \psi$, which implies  $\Delta_{\phi} \vdash_{(n,k)} (\psi \ra \theta) \ra \theta$, because 
$\vdash_{(n,k)} \psi \ra ((\psi \ra \theta) \ra \theta)$. Hence, $\Delta_{\phi} \vdash_{(n,k)} (\psi \ra \theta) \ra \theta^{\ast}$, by 
Prop. \ref{lista-formulas-utiles-2}.a). Now, 
considering that 
$\Delta_{\phi}^v \vdash_{(n,k)} (\psi \rar \theta)^{\ast}$ and $\Delta_{\phi}^v\vdash_{(n,k)} (\theta^{\ast})^{\circ}$ (because $Ax_3$ and Prop. \ref{simil-idempotencia}, resp.), it is valid that $\Delta_{\phi} \vdash_{(n,k)} \n (\theta^{\ast}) \ra \n (\psi \ra \theta)$, by 
Prop. \ref{formulas-derivadas-inpk}.d).
In addition, reasoning as in Subcase 3.2) (w.r.t $\theta$), $\Delta_{\phi}\vdash_{(n,k)} \n (\theta)^{\ast}$. Thus, 
$\Delta_{\phi}\vdash_{(n,k)} \n(\psi \rar \theta)$ (= $\phi^{v}$), as it is desired.
%

\noindent  3.5): $v(\psi)$ = $T_i$ ($0 \leq i \leq k$), $v(\theta)$ = $F_0$. 
Adapting $(\star \star)$ of Subcase 3.1) to $\theta$ it can be obtained $\Delta_{\phi}^v \vdash_{(n,k)} \theta^{\circ}$. 
Besides that, by (I.H), $\Delta_{\phi}^v \vdash_{(n,k)} \psi$ and $\Delta _{\phi}^v \vdash_{(n,k)}
\n \theta$. Considering Prop. \ref{lista-formulas-utiles-2}.g) now, it holds $\Delta_{\phi}^v \vdash_{(n,k)} \n (\psi \ra \theta)$ = $\phi^v$.

\noindent The analysis of this last subcase finishes the proof. \hfill $\Box$

\begin{lema}\label{tecnico-2}
\rm{Let $\Delta \cup \{ \psi, \theta \}$ be a subset of $L(C)$. If the following $n+k+4$ syntactic consequences are valid:

\noindent $\begin{array}{ll} {1)} &\Delta, \n(\psi^{\ast}), (\n \psi)^{\ast}  \vdash_{(n,k)} \theta \\
{ 2)}& \Delta, \n(\psi^{\ast}), \n((\n \psi)^{\ast}), (\n^2 \psi)^{\ast}\vdash_{(n,k)} \theta\\
\vdots & \vdots \\
{ n-1)} &\Delta, \n(\psi^{\ast}), \dots, \n((\n^{n-2}\psi)^{\ast}), (\n^{n-1}\psi)^{\ast}\vdash_{(n,k)} \theta\\
{ n)} &\Delta,  \n(\psi^{\ast}), \dots, \n((\n^{n-1}\psi)^{\ast}), (\n^{n}\psi)^{\ast} \vdash_{(n,k))} \theta\\
{ n+1)} &\Delta, \, \n \psi \wedge \psi, \, 
(\n \psi)^{\circ} \vdash_{(n,k)} \theta\\
{ n+2)} &   \Delta, \,  \n \psi \wedge \psi, \,  \n^2 \psi \wedge \n
\psi, \, 
(\n ^2 \psi)^{\circ} \vdash_{(n,k)} \theta\\
\vdots & \vdots \\
{ n+k-1)} &\Delta, \, \n \psi \wedge \psi,\,  \dots, \,  \n^{k-1} \psi
\wedge \n^{k-2}\psi, \,  (\n^{k-1} \psi )^{\circ} \vdash_{(n,k)} \theta\\
{ n+k)} &  \Delta, \n \psi \wedge \psi,\,  \dots, \n^{k} \psi
\wedge \n^{k-1} \psi, \, 
(\n^{k} \psi)^{\circ} \vdash_{(n,k)} \theta\\
{ n+k+1)} &  \Delta, \sim \psi, \psi^{\ast} \vdash_{(n,k)} \theta\\
{ n+k+2)} &  \Delta, \copyright \psi, \psi^{\circ} \vdash_{(n,k)} \theta\\
{ n+k+3)} &  \vdash_{(n,k)} \theta^{\ast}\\
{ n+k+4)} &  \vdash_{(n,k)} \theta^{\circ}\\
\end {array}$

\noindent Then it is valid that $\Delta \vdash_{(n,k)} \theta$.} 
\end{lema}

\dm First,
by Hypothesis ${1)}$ to $n)$ can be obtained ${ \Delta, \n (\psi^{\ast}) \vdash_{(n,k)} \theta }$ $(\star)$. Indeed:
by $n-1)$,
$\Delta, \n(\psi^{\ast}), \dots, \n((\n^{n-2}\psi)^{\ast})\vdash_{(n,k)}
(\n^{n-1}\psi)^{\ast}\rar \theta$.
Besides that, it holds that $\vdash_{(n,k)}((\n^{n-1} \psi)^{\ast})^{\ast}$ (by  Prop. \ref{simil-idempotencia}), and $\vdash \theta^{\circ}$ (by Hypothesis n+k+4)). 
 Applying all this to 
Prop. \ref{formulas-derivadas-inpk}.{d)} it is verified:

\noindent $\Delta, \n ({\psi}^{\ast}),\dots,\n((\n^{n-2}
\psi)^{\ast}) \vdash_{(n,k)}
\n \theta \ra \n ((\n^{n-1} \psi)^{\ast})$ $(\dag)$. 

\noindent It is also valid 
$\Delta,  \n(\psi^{\ast}), \dots, \n((\n^{n-2}\psi)^{\ast}) \vdash_{(n,k))} 
\n((\n^{n-1}\psi)^{\ast}) \rar \theta$, because
$Ax_{5}$, $n)$  and DT.
In addition, $\vdash_{(n,k)}
(\n((\n^{n-1}\psi)^{\ast}))^{\ast}$ , because
Proposition \ref{simil-idempotencia}  and $Ax_{11}$. Thus,
considering Prop. \ref{formulas-derivadas-inpk}.{d)} and Hyp. n+k+4) again, it holds:

\noindent $\Delta, \n ({\psi}^{\ast}),\dots,\n((\n^{n-2}
\psi)^{\ast}) \vdash_{(n,k)}
\n \theta \ra \n^2 ((\n^{n-1} \psi)^{\ast})$ $(\dag \dag)$. 

\noindent Hence, from $(\dag)$ and $(\dag \dag)$ and Corollary \ref{ejemplos-sim-n-y-o}.d):

\noindent $\Delta, \n ({\psi}^{\ast}),\dots,\n ((\n^{n-2}
\psi)^{\ast}) \vdash \n \theta \ra [ \n^2 ((\n^{n-1} \psi)^{\ast})
\wedge \n
((\n^{n-1} \psi)^{\ast})]$ $(\Diamond)$.

\noindent On the other hand, it holds $\vdash_{(n,k)} (\n \theta)^{\ast}$, because Hyp. $n + k + 3)$ and $Ax_{11}$. And, of course,
$\vdash_{(n,k)} [ \n^2 ((\n^{n-1} \psi)^{\ast})
\wedge \n
((\n^{n-1} \psi)^{\ast})]^{\circ}$. So, by Proposition \ref{formulas-derivadas-inpk}.{d)}:

\noindent $\Delta, \n({\psi}^{\ast}),\dots,\n((\n^{n-2}
\psi)^{\ast})\vdash_{(n,k)} \n[\n^2((\n^{n-1} \psi)^{\ast}) \wedge \n
((\n^{n-1} \psi)^{\ast})] \ra
\n \n \theta$. That is, 
$\Delta, \n({\psi}^{\ast}),\dots,\n((\n^{n-2}
\psi)^{\ast})\vdash_{(n,k)} (\n((\n^{n-1} \psi)^{\ast}))^{\circ} \ra
\n \n \theta$.
Thus, from Prop. \ref{lista-formulas-utiles-2}.i),
$\Delta, \n({\psi}^{\ast}),\dots,\n((\n^{n-2}
\psi)^{\ast})\vdash_{(n,k)} 
\n \n \theta$. Thus (by Hyp. n+k+3) and $Ax_9$),
$\Delta, \n({\psi}^{\ast}),\dots,\n((\n^{n-2}
\psi)^{\ast})\vdash_{(n,k)} 
\theta$ $(\Diamond \Diamond)$.

\noindent The procedure used above to prove $(\Diamond \Diamond)$ can be applied using (in decreasing order) the Hypotheses $1)$ ... $n-1)$, proving $(\star)$ (note that the formula $\n({\psi}^{\ast})$ cannot be ``suppressed'' yet). 

\noindent From $(\star)$ (and monotonicity), $\Delta, \sim \psi \vdash \n (\psi^{\ast}) \rar \theta$. Moreover, from Hyp. $n+k+1)$, it holds
$\Delta, \sim \psi \vdash \psi^{\ast}\rar \theta$. From these facts and Corollary \ref{ejemplos-sim-n-y-o}.c), it is valid
$\Delta, \sim \psi \vdash_{(n,k)}   \n (\psi^{\ast}) \vee \psi^{\ast} \rar \theta$. Now realizing that 
$\vdash_{(n,k)} \n (\psi^{\ast}) \vee \psi^{\ast}$ (because Prop. \ref{simil-idempotencia}), it is obtained
$\Delta, \sim \psi \vdash_{(n,k)}\theta$. ${\bf (I)}$.

\noindent On the other hand, 
from $n+1)$ to $n+k)$ it is valid 
$\Delta,\psi \wedge \n \psi \vdash_{(n,k)} \theta $ $(\star \star)$. The reasoning is as follows:
using ${n+k)}$ and $Ax_6$, it holds:

\noindent $\Delta, \n \psi \wedge \psi,\dots, \, \n^{k-1}
\psi \wedge \n^{k-2} \psi \vdash_{(n,k)} \n^{k} \psi \wedge
\n^{k-1} \psi \, \rar  \theta$. It is also valid 

\noindent $\Delta, \n \psi \wedge \psi,\, \dots, \n^{k-1} \psi
\wedge \n^{k-2}\psi \, \vdash_{(n,k)} (\n^{k-1} \psi )^{\circ} \rar \theta$, because Hyp. 
$n+k-1)$. Hence, 
by Corollary \ref{ejemplos-sim-n-y-o}.c) and recalling Definition \ref{conectivos-secundarios}:
 
\noindent $\Delta,
\dots,\, \n^{k-1} \psi
\wedge \n^{k-2}\psi \vdash_{(n,k)} \, 
\n(\n^{k} \psi \wedge
\n^{k-1} \psi)
\vee
(\n^{k} \psi \wedge
\n^{k-1} \psi) \,
\rar \theta$. That is,
$\Delta,
\dots,\, \n^{k-1} \psi
\wedge \n^{k-2}\psi \, \vdash_{(n,k)} (\n^{k} \psi \wedge
\n^{k-1} \psi)^{\ast}\rar \theta$.
In addition, it holds $\vdash_{(n,k)}(\n^{k-1} \psi \wedge
\n^k \psi)^{\ast}$, by Definition \ref{conectivos-secundarios}, $Ax_3$ and $Ax_{11}$. From these two facts, it holds
$\Delta,
\dots,\, \n^{k-2} \psi
\wedge \n^{k-1}\psi \, \vdash_{(n,k)}  \theta$. $(\Diamond \Diamond \Diamond)$

\noindent Adapting the reasoning applied in $(\Diamond \Diamond \Diamond)$ to the Hypotheses $n+k-2),\dots, n+1)$ (in a decreasing order, as before), it is obtained $(\star \star)$, as desired. 

\noindent From $(\star \star)$ and monotonicity it is valid $\Delta, \copyright \psi \vdash_{(n,k)}\n \psi \wedge \psi  \rar \theta $.
So (by Hyp. n+k+4), Prop. \ref{simil-idempotencia}, $Ax_{11}$ and Prop. \ref{formulas-derivadas-inpk}.{d)}), $\Delta, \copyright \psi \vdash_{(n,k)}\n \theta \rar \psi^{\circ}$.
 On the other hand, by Hyp. $n+k+2)$, it holds
$\Delta, \copyright \psi \vdash_{(n,k)}\psi^{\circ}\rar \theta$. So, 
$\Delta, \copyright \psi \vdash_{(n,k)}\n \theta \rar \n (\psi^{\circ})$ (again, by  Hyp. n+k+4), Prop. \ref{simil-idempotencia}, $Ax_{11}$ and Prop. \ref{formulas-derivadas-inpk}.{d)}). 
Thus, $\Delta, \copyright \psi \vdash_{(n,k)}\n \theta \rar ( \n (\psi^{\circ}) \wedge \psi^{\circ})$, by Corollary \ref{ejemplos-sim-n-y-o}.d). 
Therefore, 
$\Delta, \copyright \psi \vdash_{(n,k)}\n (\n (\psi^{\circ}) \wedge 
\psi^{\circ}
) \rar \n \n \theta$ (because Hyp. n+k+3), $Ax_{11}$ and Prop. \ref{formulas-derivadas-inpk}.{d)}) . That is,
$\Delta, \copyright \psi \vdash_{(n,k)}(\psi^{\circ})^{\circ}\rar \n \n \theta$. Hence, $\Delta, \copyright \psi \vdash_{(n,k)}\n \n \theta$, because Prop. 
\ref{lista-formulas-utiles-2}.c). Now, taking into account Hyp. n+k+3) and $Ax_9$, it is valid $\Delta, \copyright \psi \vdash_{(n,k)}\theta$. ${\bf (II)}$.

\noindent From ${\bf (I)}$, ${\bf (II)}$ and Corollary \ref{ejemplos-sim-n-y-o}.c), is verified
$\Delta \vdash_{(n,k)}(\copyright \psi)^{\ast} \rar \theta$. Hence, it is valid
 $\Delta \vdash_{(n,k)}\theta$, by Prop \ref{simil-idempotencia}. \hfill $\Box$

\

\noindent Thus, using Lemmas \ref{tecnico-1} and \ref{tecnico-2} it is possible to demonstrate (weak) completeness as the following result shows:

\begin{teor}\label{inpk-completitud-debil}\rm{[Weak Completeness] $\models_{(n,k)}\phi$ implies $\vdash_{(n,k)}\phi$.}
\end{teor}

\dm Suppose $\models_{(n,k)} \phi$, with $\phi$ = $\phi[\alpha_1, \dots, \alpha_m]$, and consider the set $VAL_{\phi}:=$ $\{v_t\}_{1 \leq t \leq (n+k+2)^m}$ (the set of all the $I^n P^k$-valuations effectively used to evaluate $\phi$). 
Define in $VAL_{\phi}$ the equivalence relation $\equiv_1$, as follows: for every 
$v_{t_1}, v_{t_2}\in VAL_{\phi}$,
$v_{t_1} \equiv_1 v_{t_2}$ iff, for every $\alpha_p$ {\it with $2 \leq p \leq m$}, $v_{t_1}(\alpha_p)$ = $v_{t_2}(\alpha_p)$. This relation has $(n+k+2)^{m-1}$ equivalence classes (indicated, in a general way, by $||v||$). 
Besides that, taking into account Definition \ref{conjuntos-asociados-inpk}, it holds that (given a fixed equivalence class $||v||$) $Q_p^{v_{t_1}}$ = 
$Q_p^{v_{t_2}}$, for every 
$2 \leq p \leq m$, for every pair $v_{t_1}$, $v_{t_2} \in ||v||$. 
 This allows to define the set $\Delta_1^{||v||}$:= $Q_2^{v_{t}} \cup \dots \cup Q_{m}^{v_{t}}$, being $v_t$ any element of $||v||$. 
In addition, note that every class $||v||$ has exactly $(n+k+2)$ valuations and verifies that,
for every $v_{t_1}$, $v_{t_2} \in ||v||$, $v_{t_1} \neq v_{t_2}$ implies $v_{t_1}(\alpha_1) \neq v_{t_2}(\alpha_1)$. Finally, note that, since $\models_{(n,k)}\phi$, for every $v \in VAL_{\phi}$, $\phi^v$ = $\phi$.
All these facts (together with Lemma \ref{tecnico-1}) imply that (for every $||v||$) 
the following formal proofs can be built:

\

\noindent $\begin{array}{ll} {1)} &\Delta_1^{||v||}, \n({\alpha}_1^{\ast}), (\n {\alpha}_1)^{\ast}  \vdash_{(n,k)} \phi \\
{ 1.2)}& \Delta_1^{||v||}, \n({\alpha}_1^{\ast}), \n((\n {\alpha}_1)^{\ast}), (\n^2 {\alpha}_1)^{\ast}\vdash_{(n,k)} \phi\\
\vdots & \vdots \\
{ 1.n)} &\Delta_1^{||v||},  \n({\alpha}_1^{\ast}), \dots, \n((\n^{n-1}{\alpha}_1)^{\ast}), (\n^{n}{\alpha}_1)^{\ast} \vdash_{(n,k))} \phi\\
{ 1.(n+1))} &\Delta_1^{||v||}, \, {\alpha}_1 \wedge \n {\alpha}_1,\,
(\n {\alpha}_1)^{\circ} \vdash_{(n,k)} \phi\\
{ 1.n+2))} &   \Delta_1^{||v||},\, {\alpha}_1 \wedge \n {\alpha}_1,\, \n^2 {\alpha}_1 \wedge \n {\alpha}_1, \,
(\n ^2 {\alpha}_1)^{\circ} \vdash_{(n,k)} \phi\\
\vdots & \vdots \\
{ 1.(n+k))} & \Delta_1^{||v||},\, \n {\alpha}_1 \wedge {\alpha}_1,\dots,\, \n^{k} {\alpha}_1
\wedge \n^{k-1} {\alpha}_1,\,
(\n^{k} {\alpha}_1)^{\circ} \vdash_{(n,k)} \phi\\
{ 1.(n+k+1))} &  \Delta_1^{||v||}, \sim {\alpha}_1, {\alpha}_1^{\ast} \vdash_{(n,k)} \phi\\
{ 1.(n+k+2))} &  \Delta_1^{||v||}, \copyright {\alpha}_1, {\alpha}_1^{\circ} \vdash_{(n,k)} \phi\\
\end {array}$

\noindent Moreover, by Proposition \ref{dem-bien-comportadas}, it is valid:

\noindent $\begin{array}{ll}
{ 1.(n+k+3))} & \Delta_1^{||v||} \vdash_{(n,k)} \phi^{\ast}\\
{ 1.(n+k+4))} &  \Delta_1^{||v||} \vdash_{(n,k)} \phi^{\circ}\\
\end {array}$

\noindent All the previous facts allow to apply Lemma \ref{tecnico-2} in such a way that for every $||v||$ it holds $\Delta_1^{||v||} \vdash_{(n,k)} \phi$ (there are $(n+k+2)^{m-1}$ formal proof of this type). That is, it is possible ``to eliminate'' any referrence to formulas of the form $\alpha_1^v$ in every formal proof obtained,  by means of an adequate subdivision of the set $VAL_{\phi}$, and by the application of Lemma \ref{tecnico-2}. Note here that this process can be applied one more time, reagrouping the formal proofs already obtained. So, by a new application of Lemma \ref{tecnico-2} and of Proposition \ref{dem-bien-comportadas}, any referrence to formulas of the form $\alpha_2^{v}$ can be supressed. 
The same prodedure can be applied by a finite number of times, until obtaining the following formal proofs:

\noindent $\begin{array}{ll} 
{m.1)} & \n({\alpha}_m^{\ast}), (\n {\alpha}_m)^{\ast}  \vdash_{(n,k)} \phi \\
{ m.2)}&  \n({\alpha}_m^{\ast}), \n((\n {\alpha}_m)^{\ast}), (\n^2 {\alpha}_m)^{\ast}\vdash_{(n,k)} \phi\\
\vdots & \vdots \\
{ m.n)} &  \n({\alpha}_m^{\ast}), \dots, \n((\n^{n-1}{\alpha}_m)^{\ast}), (\n^{n}{\alpha}_m)^{\ast} \vdash_{(n,k))} \phi\\
{ m.(n+1))} & \n {\alpha}_m \wedge {\alpha}_m,\,
(\n {\alpha}_m)^{\circ} \vdash_{(n,k)} \phi\\
{ m.(n+2))} &\n {\alpha}_m \wedge {\alpha}_m,\, \n^2 {\alpha}_m \wedge \n {\alpha}_m,\,
(\n ^2 {\alpha}_m)^{\circ} \vdash_{(n,k)} \phi\\
\vdots & \vdots \\
{ m.(n+k))} &  \n {\alpha}_m \wedge {\alpha}_m,\, \dots, \n^{k} {\alpha}_m
\wedge \n^{k-1} {\alpha}_m,\, 
(\n^{k} {\alpha}_m)^{\circ} \vdash_{(n,k)} \phi\\
{ m.(n+k+1))} & \sim {\alpha}_m, {\alpha}_m^{\ast} \vdash_{(n,k)} \phi\\
{ m.(n+k+2))} &  \copyright {\alpha}_m, {\alpha}_m^{\circ} \vdash_{(n,k)} \phi\\
{ m.(n+k+3))} &  \vdash_{(n,k)} \phi^{\ast}\\
{ m.(n+k+4))} &  \vdash_{(n,k)} \phi^{\circ}\\
\end {array}$

\noindent Applying Lemma \ref{tecnico-2} and Proposition \ref{dem-bien-comportadas} for a last time, $\vdash_{(n,k)}\phi$. \hfill $\Box$

\

\noindent Note that, in the proof developed above, {\it all the ${(n+k+2)}^m$ valuations of $VAL_{\phi}$} are needed to obtain the formal proofs that allow to demonstrate 
$\vdash_{(n,k)}\phi$.

\

\noindent Theorems \ref{inpk-correccion-debil} and \ref{inpk-completitud-debil} prove {\it weak adequacity}: $\models_{(n,k)}\phi$ iff $\vdash_{(n,k)}\phi$. This result can be improved:

\begin{teor}\label{inpk-completitud-fuerte}\rm{[Strong Adequacity]: for every $\Gamma \cup \{\phi\} \su L(C)$, $\Gamma \models_{(n,k)}\phi$ iff $\Gamma \vdash_{(n,k)}\phi$.}
\end{teor}
\dm  By Proposition \ref{models-finitaria-dt}, $\models_{(n,k)}$ verifies Semantics Deduction Theorem and is finitary. Moreover, by the definition of formal proof used in this paper, $\vdash_{(n,k)}$ is finitary,  and (by Theorem \ref{teor-deduccion-sintactico}) it verifies Sintactic Deduction Theorem, as was already mentioned. From all this facts, and taking into account that both $\models_{(n,k)}$ and $\vdash_{(n,k)}$ are monotonic, strong adequacity is demonstrated. \hfill $\Box$ 

\section{Concluding remarks}
Despite its interest as a general result (for a countable, non-lineal family of logics), the adequate axiomatics shown here  can be applied in different ways. First of all, a natural problem to be solved is the {\it independence} of the axiomatics presented here and it is part of a future work. 

\

\noindent On the other hand, 
another of the possible uses of this axiomatics is the study of {\it algebraizability} of the $I^n P^k$-logics. It is worth to comment here that $I^1 P^0$ is algebraizable (see \cite{set:car:95}), as in the case of $I^0 P^1$ (this fact was already indicated). Moreover, in \cite{fer:05} it was demonstrated that {\it all the logics of ${\mathbb{I}}^n {\mathbb{P}}^k$} are algebraizable. So, the properties of the class of algebras associated to each $I^n P^k$-logic deserve to be investigated. By the way, the class of algebras associated to $I^0 P^1$ was already studied in \cite{lew:mik:sch:94} and in \cite{pyn:95}. In both works, the axiomatics obtained for this logic are very useful for the study of the so-called {\it class of $P^1$-algebras}. This is because there is a connection between the axiomatics of an algebraizable logic and its equivalent algebraic semantics, cf. \cite{blo:pig:89}. As a generalization of this fact, the axiomatics shown here would allow to study the different classes  of (say) $I^n P^k$-algebras in a more efficient way. 

\

\noindent Finally, note this fact about the complexity of the formulas: given a fixed logic $I^n P^k$, every formula $\phi \in L(C)$ with complexity $Comp(\phi) \geq max\{n,k\}$ behaves ``in a classical way'' (this fact is related to the inclusion of $Ax_5$ and $Ax_6$ in the axiomatics presented in this paper). This would suggest to define a special kind of logics: the family $\mathbb{SC}$ of ``stationary classically logics''. Obviously, ${\mathbb{I}}^n {\mathbb{P}}^k$ would be a particular subclass of $\mathbb{SC}$. The study of the latter class deserves special attention in a future research.

\end{document}